\documentclass[twoside]{article}

\usepackage{amsfonts}
\usepackage{amsthm}
\usepackage{amsmath}
\usepackage{amssymb}

\title{\bf On some properties of travelling water waves with vorticity}
\author{Eugen Varvaruca}\date{}

\newtheorem{theorem}{Theorem}[section]

\numberwithin{equation}{section}

\newcommand{\Om}{\Omega}
\newcommand{\ovo}{\overline\Om}

\newcommand{\om}{\omega}

\newcommand{\be}{\begin{equation}}
\newcommand{\ee}{\end{equation}}

\newcommand{\mcs}{{\mathcal S}}

\newcommand{\mcb}{{\mathcal B}}

\newcommand{\hg}{\hat\Gamma}

\newcommand{\psix}{\psi_x}
\newcommand{\psiy}{\psi_y}

\newcommand{\bdr}{\mathbb{R}}

\newcommand{\bese}{\begin{subequations}}
\newcommand{\ese}{\end{subequations}}

\begin{document}
\maketitle \centerline{Department of Mathematical Sciences,
University of Bath} \centerline{Claverton Down, Bath BA2 7AY, United
Kingdom}\centerline{Email address: {\tt mapev@maths.bath.ac.uk}}
\centerline{}

\begin{abstract} We prove that for a large class of vorticity
functions the crest of a corresponding travelling water wave is
necessarily a point of maximal horizontal velocity. We also show
that for waves with nonpositive vorticity the pressure in the flow
is everywhere larger than the atmospheric pressure. A related a
priori estimate for waves with nonnegative vorticity is also given.
\end{abstract}

\section{Introduction}
In this article we consider the classical hydrodynamical problem
concerning travelling two-dimensional gravity water waves with
vorticity. This problem has attracted considerable interest in
recent years, starting with the systematic study of Constantin and
Strauss \cite{CoS} on periodic waves of finite depth.

 The problem arises from the following
physical situation. A wave of permanent form moves with constant
speed on the surface of an incompressible flow, the bottom of the
flow domain being horizontal. With respect to a frame of reference
moving with the speed of the wave, the flow is steady and occupies a
fixed region $\Om$ in the upper half of the $(x,y)$-plane, which
lies between the real axis $\mcb:=\{(x,0):x\in\bdr\}$ and some a
priori unknown free surface
$\mathcal{S}:=\{(x,\eta(x)):x\in\mathbb{R}\}$, where $\eta$ is a
periodic function.
 Since the fluid is incompressible, the
flow can be described by a (relative) stream function $\psi$
 which is periodic in the horizontal direction and
  satisfies the following equations and boundary conditions:
\begin{subequations}\label{pb}
\begin{align}
& \Delta\psi=-\gamma(\psi)\quad\text{in }\Om,\label{p0}
\\& \psi=B\quad \text{on }\mcb,\label{p1}\\&\psi=0\quad
\text{on }{\mathcal S}, \label{p2}\\&\vert\nabla\psi\vert^{2}+2gy=Q
\quad\text{on } {\mathcal S}\label{p3},\\
&\psi_y<0 \quad\text{in }\Om,\label{p4}
 \end{align}
\end{subequations}
where $B,\,g$ and $Q$ are positive constants. The equation
(\ref{p0}) involves a vorticity function $\gamma:[0,B]\to\bdr$ and
expresses the fact that the vorticity of the flow $\om:=-\Delta\psi$
and the stream function $\psi$ are functionally dependent. The
equations (\ref{p1}) and (\ref{p2}) mean that the bottom and the
free surface are streamlines, while (\ref{p2}) means that the
pressure at the surface of the flow is a constant.  The relative
velocity of the fluid particles is given by $(\psi_y, -\psi_x)$, and
the requirement (\ref{p4}) is motivated by physical considerations.
It is customary \cite{CoS} to assume that the constants $g, B$ and
the vorticity function $\gamma$ are given. The problem consists in
determining the curves $\mcs$ for which there exists a function
$\psi$ in $\Om$ which satisfies (\ref{p0})-(\ref{p4}) for some value
of the parameter $Q$. For a full justification of the equivalence
between the problem of seeking solution triples $(\mcs, \psi, Q)$ of
(\ref{pb}) and that of seeking travelling-wave solutions of the
two-dimensional Euler equations, the reader is referred to
\cite{CoS}.

 When $\gamma\equiv 0$, the corresponding flow is called
irrotational. Nowadays the mathematical theory dealing with this
situation contains a wealth of results, mostly obtained during the
last three decades, concerning the existence of large amplitude
solutions and their properties. Global bifurcation theories were
given for various types of waves (periodic or solitary of finite
depth, periodic of infinite depth) by Keady and Norbury \cite{KN},
and by Amick and Toland \cite{AT1,AT2}. Moreover, it was shown by
Toland \cite{T78} and McLeod \cite{ML} that in the closure of these
continua of solutions there exist waves with stagnation points at
their crests, a stagnation point being one at which the relative
fluid velocity is zero, i.e.\ $|\nabla\psi|=0$. The existence of
such waves, called `extreme waves', was predicted by Stokes
\cite{S}, who also conjectured that their profiles necessarily have
 corners with included angle of $120^\circ$ at the crests. This
conjecture was proved independently by Amick, Fraenkel and Toland
\cite{AFT} and by Plotnikov \cite{P82}. Recently, the method of
\cite{AFT} was simplified and extended in \cite{EV}.

On the other hand, when $\gamma\not\equiv 0$, the flow is called
rotational or with vorticity, and significant advances in the
corresponding mathematical theory have been made only in the last
few years. The existence of global continua of smooth solutions was
proved by Constantin and Strauss \cite{CoS} for the periodic finite
depth problem, and by Hur \cite{H} for the related problem of
periodic waves of infinite depth. For the solutions found in
\cite{CoS,H} the wave profiles are monotone between each crest and
trough and have a vertical axis of symmetry.
 Of particular significance is the fact that
  the continuum of solutions found in \cite{CoS} contains waves for
  which the value of $\max_{\overline\Om}\psi_y$ is arbitrarily close to
  $0$. Thus it is natural to expect that, as in the irrotational case,
waves with stagnation points, referred to as `extreme waves' exist
for many vorticity functions, and that they can be obtained as
limits, in a suitable sense, of certain sequences of smooth
solutions found in \cite{CoS}. In the case of constant vorticity,
numerical evidence \cite{KS, SS, TSP, VB1, VB2, VB3} strongly points
to the existence of extreme waves for any negative vorticity and for
small positive vorticity, and also indicates that, for large
positive vorticity, continua of solutions bifurcating from the
trivial solution develop into overhanging profiles (a situation
which is not possible in the irrotational case, see \cite{EV3} for
references) and do not approach extreme waves.

One of the questions addressed in this article concerns the location
of the points at which the maximum over $\ovo$ of the relative
horizontal velocity $\psiy$ is attained for smooth waves with
vorticity. In the irrotational case, the crest of the wave is
necessarily such a point, see Toland \cite{T}. Very recently,
Constantin and Strauss \cite[Theorem 4.1]{CoS1} showed that this is
also the case for the waves in the continuum found in \cite{CoS}
under the conditions that $\gamma$ is a nonpositive constant which
satisfies a smallness condition involving $B$ and $g$. Here we
substantially extend this result, in that we only require that the
vorticity function $\gamma$ satisfies $\gamma\leq 0$ and
$\gamma'\geq 0$ everywhere on $[0,B]$, while our assumptions on the
wave profiles are also more general. An immediate consequence of
this result is that, whenever $\gamma\leq 0$ and $\gamma'\geq 0$,
the continuum of solutions in \cite{CoS} contains waves for which
the values of $|\nabla\psi|$ at their crests are arbitrarily close
to $0$. Thus in this case the existence of waves with stagnation
points at their crests is to be particularly expected.

Another contribution of this article is that we establish some new a
priori bounds for waves corresponding to vorticity functions
$\gamma$ which do not change sign, without any assumptions on
$\gamma'$. When $\gamma\leq 0$, the estimate in question means that
the pressure everywhere in the fluid is larger than the atmospheric
pressure; this estimate is the essential ingredient in the proof of
the previously mentioned result concerning the location of the
points where $\max_{\overline\Om}\psi_y$ is attained. When
$\gamma\geq 0$, a slightly different, but related, estimate is
given. Both these estimates play an essential role in the
investigation in \cite{EVv} concerning the existence of extreme
waves with vorticity and the Stokes conjecture.

Analogous results to those of this article hold in the case of
periodic rotational waves of infinite depth. They will be presented,
together with some applications, in a subsequent article.

\section{The Main Results}

We always deal with classical solutions of (\ref{pb}), in the sense
that $\gamma\in C^1([0,B])$, $\eta\in C^3(\bdr)$, $\psi\in
C^3(\overline\Om)$. We assume that $\eta$ is a periodic function of
minimal period $2L$, and that $\psi$ is $2L$-periodic in the
horizontal direction. However, we do not assume that $\eta$ has
exactly one local maximum and one local minimum per minimal period.

Let $\hg :[0,B]\to\bdr$ be given by \be
\hg(s)=\int_0^s\gamma(t)\,dt\quad\text{for all }s\in[0,B].\ee (Note
that in \cite{CoS} a function $\Gamma$ is considered which is
related to $\hg$ by $\hg(s)=-\Gamma(-s)$. The quantity of interest
both here and there is $\hg(\psi)$, which is denoted there by
$-\Gamma(-\psi)$; we find our notation more convenient.) Let us also
consider the function $R:\overline\Om\to\bdr$ given by \be
R=\frac{1}{2}\vert\nabla\psi\vert^{2}+gy-\frac{1}{2}Q+\hg(\psi).\label{r}\ee
The function $R$ is (up to a constant) the negative of the pressure
in the fluid domain, see \cite{CoS}.

Our next result shows that when $\gamma$ is everywhere nonpositive
the pressure in the fluid domain is larger than the atmospheric
pressure. The same conclusion was obtained in \cite[Example
3.1]{CoS1} under much more restrictive assumptions on $\gamma$ and
$\psi$, namely that
\[\gamma\leq 0,\quad\gamma'\leq 0\quad\text{ and }\quad\text{$-\psiy(x,0)\gamma(B)\geq -g$
for all $x\in\bdr$.}\]

\begin{theorem}\label{tmr}
Suppose that $\gamma(s)\leq 0$ for all $s\in [0,B]$. Then $R\leq 0$
in $\overline\Om$.
\end{theorem}

The importance of the inequality $R\leq 0$ in $\overline\Om$ in
relation to the monotonicity of $\psiy$ along the free surface
$\mcs$ was first recognized for waves with vorticity by Constantin
and Strauss \cite[Proposition 3.4]{CoS1}. We give here a slightly
more general statement of their result and a somewhat more direct
proof.

\begin{theorem}\label{maxs} Let $\eta:\bdr\to\bdr$ be such that that there exists
$N\in\mathbb{N}$ and points $x_0<x_1<\dots<x_{2N}=x_0+2L$ with the
property that $\eta'(x_j)=0$ for all $j\in\{0,\dots,2N\}$, $\eta$ is
strictly increasing on $[x_{2j},x_{2j+1}]$ for all $j\in\{0,\dots,
N-1\}$ and $\eta$ is strictly decreasing on $[x_{2j-1},x_{2j}]$ for
all $j\in\{1,\dots,N\}$. Suppose that $R\leq 0$ in $\ovo$. Then the
function $x\mapsto \psiy(x,\eta(x))$ is increasing on
$[x_{2j},x_{2j+1}]$ for all $j\in\{0,\dots, N-1\}$ and decreasing on
$[x_{2j-1},x_{2j}]$ for all $j\in\{1,\dots,N\}$. Therefore,
$\max_\mcs \psiy$ is attained at the points of maximal height on
$\mcs$.
\end{theorem}

The preceding result leads with little effort to one concerning the
location of the points where $\max_{\ovo}\psiy$ is attained. It
significantly improves a result of Constantin and Strauss
\cite[Theorem 4.1]{CoS1}, where the same conclusion was obtained,
for a more restrictive class of wave profiles, under the assumption
that $\gamma$ is a nonpositive constant which satisfies
\[g^2\geq 2g(-2B\gamma^3)^{1/2}-2B\gamma^3.\]

\begin{theorem}\label{mal} Let $\eta:\bdr\to\bdr$ be as in Theorem \ref{maxs}.
Suppose that $\gamma(s)\leq 0$ and $\gamma'(s)\geq 0$ for all
$s\in[0,B]$. Then $\max_{\ovo}\psiy$ is attained at the points of
maximal height on $\mcs$.
\end{theorem}

The next result gives a new estimate in the case when $\gamma$ is
everywhere nonnegative, which is in the same spirit as that of
Theorem \ref{tmr}. Let us consider the function $T:\ovo\to\bdr$
given by \be T:=R-\varpi\psi,\ee where $R$ is given by (\ref{r}) and
\be \varpi:=\frac{1}{2}\max_{s\in[0,B]}\gamma(s).\label{defa}\ee

\begin{theorem}\label{ppo} Suppose that $\gamma(s)\geq 0$ for all $s\in[0,B]$.
Then $T\leq 0$ in $\ovo$.
\end{theorem}

\section{Proofs of the Main Results}

 A simple calculation shows that, everywhere in $\ovo$,
\bese\label{grad}\begin{align}R_x&=\psiy\psi_{xy}-\psix\psi_{yy},
\\R_y&=\psix\psi_{xy}-\psiy\psi_{xx}+g,\end{align}\ese
and \be\qquad\Delta
R=2\psi_{xy}^2-2\psi_{xx}\psi_{yy}.\label{lapl}\ee

\begin{proof}[Proof of Theorem \ref{tmr}] Note that $R=0$ everywhere
on the free surface $\mcs$. We claim that the maximum of $R$ over
$\ovo$ must be attained on $\mcs$.

Observe first that, since $R_y=g>0$ everywhere on the bottom $\mcb$,
$\max_{\ovo} R$ cannot be attained anywhere on $\mcb$.

Suppose now that $\max_{\ovo} R$ is attained at some point $A$ in
$\Om$. Then necessarily
\[R_x(A)=0,\quad R_y(A)=0,\quad \Delta R(A)\leq 0.\]
It follows from this, (\ref{grad}) and (\ref{lapl}) that
\bese\begin{align}\psiy(A)\psi_{xy}(A)&=\psix(A)\psi_{yy}(A),\label{a}\\\psix(A)\psi_{xy}(A)&<\psiy(A)\psi_{xx}(A),\label{b}\\
\psi_{xy}^2(A)&\leq\psi_{xx}(A)\psi_{yy}(A).\label{c}
\end{align}\ese
Since (\ref{p4}) holds, it follows that $\psi_y(A)<0$. We now
distinguish two cases, depending on whether or not $\psi_{yy}(A)=0$.

 If
$\psi_{yy}(A)=0$, then (\ref{a}) implies that $\psi_{xy}(A)=0$. It
then follows from (\ref{b}) that $\psi_{xx}(A)<0$, and hence
$\gamma(\psi(A))=-\Delta\psi(A)>0$, which contradicts the assumption
that $\gamma(s)\leq 0$ for all $s\in[0,B]$.

If $\psi_{yy}(A)\neq 0$, then it follows from (\ref{a}) and
(\ref{b}) that
\[\frac{\psi_y(A)\psi_{xy}^2(A)}{\psi_{yy}(A)}<\psiy(A)\psi_{xx}(A).\]
It then follows from this and (\ref{c}) that $\psi_{yy}(A)<0$. We
now deduce from (\ref{c}) that $\psi_{xx}(A)\leq 0$, and therefore
$\gamma(\psi(A))=-\Delta\psi(A)>0$, which  again contradicts the
assumption that $\gamma(s)\leq 0$ for all $s\in[0,B]$.

We conclude that the maximum of $R$ over $\ovo$ must be attained on
$\mcs$, which implies that $R\leq 0$ in $\ovo$. This completes the
proof of Theorem \ref{tmr}.
\end{proof}

\begin{proof}[Proof of Theorem \ref{maxs}] The proof is based on a
remarkable, though straightforward to verify, identity observed by
Toland \cite{T} in the irrotational case and by Constantin and
Strauss \cite{CoS1} in the general case:
\be\frac{d}{dx}\left[\frac{1}{2}\psi_y^2(x,\eta(x))\right]=R_x(x,\eta(x))\quad\text{for
all }x\in\bdr.\label{pre}\ee Since $R\leq 0$ in $\ovo$ and $R=0$ on
$\mcs$, the required result concerning the monotonicity of $x\mapsto
\psiy(x,\eta(x))$ is immediate from (\ref{pre}). It follows that \be
\max_\mcs \psiy=\max_{j\in\{0,\dots,2N\}}\psiy(x_j, \eta(x_j)).\ee
But for every $j\in\{0,\dots,2N\}$, $\psix(x_j, \eta(x_j))=0$ and
therefore
\[\psiy(x_j, \eta(x_j))=-(Q-2g\eta(x_j))^{1/2}.\]
Hence $\max_{\ovo}\psiy$ is attained at the points of maximal height
on $\mcs$. This completes the proof of Theorem \ref{maxs}.
\end{proof}

\begin{proof}[Proof of Theorem \ref{mal}] It follows from (\ref{p0})
that
\[\Delta\psiy=-\gamma'(\psi)\psiy\quad\text{in }\Om.\]
Since $\psiy< 0$ in $\Om$ and $\gamma'(s)\geq 0$ for all
$s\in[0,B]$, it follows immediately from the Maximum Principle that
$\max_{\ovo}\psiy$ cannot be attained anywhere in $\Om$.

We now show that $\max_{\ovo}\psiy$ cannot be attained anywhere on
$\mcb$ either. This is trivial when $\gamma(B)<0$, since then
$\psi_{yy}=-\gamma(B)>0$ everywhere on $\mcb$. When $\gamma(B)=0$,
we use the Reflection Principle. Let $\tilde\gamma:[0,2B]\to\bdr$ be
an extension of $\gamma$ such that $\tilde\gamma(s)=-\gamma(2B-s)$
for all $s\in (B,2B]$. Let $\Om^R$ be the reflection of $\Om$ into
$\mcb$,
\[\widetilde\Om:=\Om\cup\mcb\cup\Om^R,\]
and $\tilde\psi:\widetilde\Om\to\bdr$ be an extension of $\psi$ such
that $\tilde\psi(x,y)=2B-\psi(x,-y)$ for all $(x,y)\in\Om^R$. Then
$\tilde\gamma\in C^1([0,2B])$, $\tilde\psi\in C^2(\widetilde\Om)$
and
\[\Delta\tilde\psi=-\tilde\gamma(\tilde\psi)\quad\text{in }\widetilde\Om.\]
Since $\tilde\psi_y<0$ in $\widetilde\Om$ and $\tilde\gamma'(s)\geq
0$ for all $s\in[0,2B]$, the Maximum Principle yields the required
result.

We conclude that $\max_{\ovo}\psiy$ is necessarily attained
somewhere on $\mcs$. Next note that, since $\gamma(s)\leq 0$ for all
$s\in[0,B]$, Theorem \ref{tmr} shows that $R\leq 0$ in $\ovo$. An
application of Theorem \ref{maxs} now yields that $\max_{\ovo}\psiy$
is attained at the points of maximal height on $\mcs$. This
completes the proof of Theorem \ref{mal}.

\end{proof}

\begin{proof}[Proof of Theorem \ref{ppo}]
 Note first that $T_y=g-\varpi\psi_y>0$ everywhere on $\mcb$, so that the maximum of $T$ over $\ovo$ cannot be
 attained anywhere on $\mcb$.

 Next note
from (\ref{lapl}) that
\[\Delta R\geq -\frac{1}{2}\gamma^2(\psi)\quad\text{in }\Om.\]
Since
\[\Delta T=\Delta R+\varpi\gamma(\psi),\]
it is immediate, upon using (\ref{defa}) and the assumption that
$\gamma(s)\geq 0$ for all $s\in[0,B]$, that $T$ is a subharmonic
function in $\Om$. Therefore the maximum of $T$ over $\overline\Om$
cannot be attained anywhere in $\Om$.

 We conclude that $\max_{\ovo}T$ must be attained somewhere on $\mcs$.
 Since $T=0$ everywhere on $\mcs$, it follows that $T\leq 0$ in
 $\ovo$. This completes the proof of Theorem \ref{ppo}.
\end{proof}

\end{document}